\documentclass[12pt]{amsart}

\usepackage{amsmath,amsthm,amssymb}
\usepackage{verbatim}
\usepackage{epsfig}
\usepackage{float}
\usepackage{mathtools}
\usepackage{algorithm}
\usepackage{algorithmic}
\usepackage{MnSymbol}
\usepackage{graphicx}
\usepackage{mathtools}
\usepackage{subcaption}
\usepackage[left=2.8cm,top=3cm,right=2.8cm,bottom=3cm]{geometry}

\newtheorem{theorem}{Theorem}

\newcommand{\G}{\mathbb{G}}
\newcommand\eps{\varepsilon}
\newcommand{\Prob}{\mathbb{P}}
\DeclarePairedDelimiter{\ceil}{\lceil}{\rceil}
\DeclarePairedDelimiter{\floor}{\lfloor}{\rfloor}

\begin{document}

\title{A survey of graph burning}\thanks{The author was supported by NSERC}

\author{Anthony Bonato}
\email{abonato@ryerson.ca}

\begin{abstract}
Graph burning is a deterministic, discrete-time process that models how influence or contagion spreads in a graph. Associated to each graph is its burning number, which is a parameter that quantifies how quickly the influence spreads. We survey results on graph burning, focusing on bounds, conjectures, and algorithms related to the burning number. We will discuss state-of-the-art results on the burning number conjecture, burning numbers of graph classes, and algorithmic complexity. We include a list of conjectures, variants, and open problems on graph burning.
\end{abstract}


\maketitle

\section{Introduction}

The spread of influence is a key topic in network science, focusing on the propagation of emotion, members, or contagion in social networks; see \cite{k}. Internet memes, for example, appear and spread quickly across social networks like Facebook, Twitter, and Instagram.  An elementary rule is that influence spreads from a vertex to each of its neighbors. While there is a source vertex from which the influence originates, other sources appear over time in various locations in the network.

Graph burning is a simplified model for the spread of influence in a network. Associated with the process is a parameter, the burning number, which quantifies the speed at which the influence spreads to ever vertex. The smaller the burning number is, the faster an influence can be spread in the network. Graph burning is defined as follows.  Given a finite, simple, undirected graph $G$, the burning process on $G$ is a discrete-time process. Vertices may be either unburned or burned throughout the process. Initially, in round $t=0$ all vertices are unburned. At each round $t \geq 1$, one new unburned vertex is chosen to burn, if such a vertex is available. We call such a chosen vertex a \emph{source}. If a vertex is
burned, then it remains in that state until the end of the process. Once a vertex is burned in round $t$, in round $t+1$ each of its unburned neighbors becomes burned. The process ends in a given round when all
vertices of $G$ are burned. We emphasize that sources are chosen in each round for which they are available.

Note that the burning process may be highly dependent on the choice of sources. For example, in a path, burning spreads more slowly from a source that is an end-vertex than from a central vertex. Hence, the strategic choice of sources is critical when minimizing the length of the process.

The burning number corresponds to an optimal choice of sources throughout the process. The \emph{burning number} of a graph $G$, denoted by $b(G)$, is the minimum number of rounds needed for the process to end. The parameter $b(G)$ is well-defined, as in a finite graph, there are only finitely many choices for the sources. The sources that are chosen over time are referred to as a \emph{burning sequence}; a shortest such sequence is called \emph{optimal}. Hence, optimal burning sequences have length $b(G).$

For an elementary example of graph burning, consider the path $P_4$ with vertices $v_i$, where $1\le i \le 4.$  In this example, the sequence $(v_2, v_4)$ is an optimal burning sequence; see Figure~\ref{CH9p4}.
\begin{figure}[h]
\begin{center}
\epsfig{figure=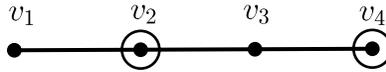}
\caption{Burning the path $P_4.$ The circled vertices are sources.} \label{CH9p4}
\end{center}
\end{figure}
We also observe in this example that optimal burning sequences may not be unique. The sequence $(v_1,v_3)$ is the other optimal burning sequence for $P_4.$

Graph burning is contained within the area of graph searching. For example, graph burning is reminiscent but distinct from the \emph{Firefighter Problem}, where a set of firefighters block burned vertices from spreading; see \cite{fin} for a survey.

Since graph burning was first introduced in \cite{BJR0,BJR,thez}, a number of results, conjectures, and algorithms have been emerged on the topic in over two dozen papers. The purpose of the present article is to survey the main results on the topic, paying attention to the central topics and questions. While we make an effort to be self-contained, as results appear in the literature, we do not claim all current topics on graph burning are represented. We view the survey as both an entry point to graph burning, and also a one-stop-shop reference for experts on the topic.

The paper is organized as followed. In Section~\ref{seccon}, we discuss one of the main open problems in the field, the burning number conjecture. While the conjecture is unresolved, we consider the best known upper bounds on the burning number. We discuss graphs families, such as spiders and caterpillars, where the conjecture is known to hold. In Section~\ref{secclass}, we consider the burning number in a variety of graph classes. We highlight the best known results on various graph products, grids, and hypercubes. We consider graph burning on binomial random graphs, generalized Petersen graphs, and theta graphs. Computational complexity results are presented in Section~\ref{seccomp}. While computing the burning number was known to be \textbf{NP}-complete early on in its formulation, a number of newer results have been discovered. We consider \textbf{NP}-completeness for graph burning in a several, restrictive graph families, and also consider research on approximation algorithms and heuristics. The article finishes with a collection of open problems, conjectures, and variants of the process.

All graphs we consider are simple, finite, and undirected, unless otherwise stated. For a vertex $v$ and a non-negative integer $k$, the {\em $r$th closed neighborhood} $N_k[v]$ of $v$ is defined as the set of all vertices within distance $k$ of $v$, including $v$ itself. In the case $k=1,$ we write $N_1[v]=N[v].$ The distance between vertices $u$ and $v$ is denoted by $d(u,v).$ If $G$ is a graph and $u$ is a vertex of $G$, then the \emph{eccentricity} of $u$ is defined as $\max\{d(u, v): v\in V(G)\}$. The \emph{radius} of $G$ is the minimum eccentricity over the set of all vertices in $G$. For background on graph
theory, see \cite{west}. For additional background on graph searching, the reader is directed to \cite{BN,bp}.

\section{Burning number conjecture}\label{seccon}

The burning number conjecture is one of the main unanswered questions on graph burning. Before we state the conjecture, we introduce alternative characterizations of graph burning in terms of neighbor sets and trees.

We first consider a characterization via a certain set equation first derived in \cite{BJR}. If $(x_1, x_2 \ldots, x_k)$ is a burning sequence for a given graph $G$, then a source at $x_i,$ where $1\le i \le k,$ will burn only all the vertices within distance
$k-i$ from $x_i$ by the end of the $k$-th step. Each vertex $v\in V(G)$ must be either a source or burned from at least one of the sources by the $k$-th
round. Further, for each pair $i$ and $j$, with $1\leq i < j\leq k$, we must have $d(x_i, x_j) \geq j-i$; otherwise, if $d(x_i, x_j) = l < j-i$, then $x_j$ will be burned at round $l + i <
j$.  Therefore, $(x_1, x_2, \ldots, x_k)$ is a burning
sequence for $G$ if and only if for each pair $i$ and $j$, with $1\leq i < j\leq k$, $d(x_i, x_j) \geq j-i$, and the following set equation holds:
$$
N_{k-1}[x_1] \cup N_{k-2}[x_2]\cup \ldots \cup N_0[x_k] = V(G).
$$

The following theorem provides another characterization of the burning number, and connects it with a prescribed covering problem by trees. The \emph{depth} of a vertex in a rooted tree is the number of edges in a shortest path from the vertex to the tree's root. The \emph{height} of a rooted tree $T$ is the greatest depth in $T$. A \emph{rooted tree
partition} of $G$ is a collection of rooted trees which are subgraphs of $G$, with the property that the vertex sets of the trees partition $V(G)$.

\begin{theorem}[\cite{BJR}]\label{decomp}
Burning a graph $G$ in $k$ steps is equivalent to finding a rooted tree partition into $k$ trees $T_1, T_2,\ldots, T_k$, with heights at most $(k-1), (k-2), \ldots, 0$, respectively such that for
every $1\le i , j \le k$ the distance between the roots of $T_i$ and $T_j$ is at least $|i-j|$.
\end{theorem}

The following theorem is a corollary of Theorem~\ref{decomp}, and is useful for determining the burning number of a graph, as it reduces the problem of burning a graph to burning its spanning trees. Note that for a spanning subgraph $H$ of
$G$, it is evident that $b(G)\leq b(H)$ (although this hereditary property does not hold in general if we consider subgraphs or induced subgraphs).

\begin{theorem}[Tree Reduction Theorem, \cite{BJR}]\label{spann}
For a graph $G$ we have that
$$b(G) = \min\{b(T): T \text{ is a spanning subtree of } G \}.$$
\end{theorem}

Paths play an important role in graph burning.

\begin{theorem}[\cite{BJR}]\label{path}
For a path $P_n$ on $n$ vertices, we have that $b(P_n) = \lceil n^{1/2}\rceil$.
\end{theorem}
Note that for any graph $G$ with radius $r$ and diameter $d$, we have that $$\lceil (d+1)^{1/2}\rceil\leq b(G)\leq r+1.$$ The bounds are tight, with the lower bound achieved by paths.

We say that a graph $G$ of order $n$ is \emph{well-burnable} if $b(G) \le \lceil n^{1/2}\rceil$. Theorem~\ref{path} tells us that paths are well-burnable, and as an immediate corollary, so is a graph with a Hamiltonian path. The following conjecture, first appearing in \cite{BJR}, states that \emph{every} graph is well-burnable.

\medskip

\noindent \textbf{Burning number conjecture}: For a connected graph $G$ of order $n$, $$b(G)\leq \lceil {n}^{1/2}\rceil.$$ \vspace{0.1in}

\medskip

If the burning number conjecture holds, then paths are examples of connected graphs with largest burning number. By the Tree Reduction Theorem, the conjecture holds if it holds for trees. Note that we require $G$ to be connected here, as otherwise, the burning number can be as large as $|V(G)|,$ as in the case for a graph with no edges.

The conjecture has resisted attempts at its resolution, although various upper bounds on the burning number are known. In \cite{bessy2}, it was proved that $$b(G)\leq \sqrt{\frac{32}{19}\cdot \frac{n}{1-\epsilon}}+\sqrt{\frac{27}{19\epsilon}}$$ and
$$b(G)\leq \sqrt{\frac{12n}{7}}+3\approx 1.309 \sqrt{n}+3$$ for every connected graph $G$ of order $n$ and every $0<\epsilon<1$. These bounds were improved in \cite{LL}, who proved the best known upper bound: $$b(G)\le \bigg \lceil
\frac{-3+\sqrt{24n+33}}{4} \bigg \rceil.$$

While the burning number conjecture is open for general graphs, it known to hold for a number of graph classes. We summarize results for such classes here. A \emph{spider} is a tree with at most one vertex of degree $3.$ A \emph{caterpillar} is a tree where deleting all vertices of degree 1 leaves a path.
\begin{figure}[h]
\begin{center}
\epsfig{figure=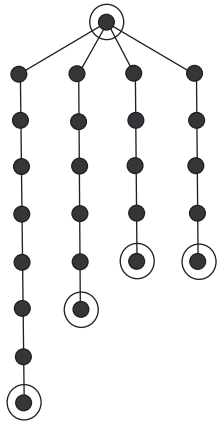}
\caption{A spider with an optimal burning sequence represented by circled vertices.} \label{CH9p4}
\end{center}
\end{figure}

As proven in \cite{bl} and independently in \cite{das}, spiders are well-burnable. As proven in \cite{liu1} and independently in \cite{hil}, caterpillars are well-burnable. For $p \ge 1$, a $p$-\emph{caterpillar} is a tree where there is a path $P$, such that each vertex is distance at most $p$ to $P$. Note that a $1$-caterpillar is a caterpillar. It was show in \cite{hil} that a $2$-caterpillar is well-burnable, although the burning number conjecture remains open for $p$-caterpillars with $p\ge 3.$ In \cite{hil}, it was shown that $p$-caterpillars with at least $2 \lceil n^{1/2} \rceil - 1$ vertices of degree one are well-burnable. In \cite{kam}, it was proven that for any graphs with minimum degree $\delta \ge 23$ are well-burnable. Although this result encompasses a large class of graphs, it omits the class of trees.

\section{Burning graph classes}\label{secclass}

The burning number has been studied in a number of graph classes, such as graph products, grids, random graphs, and certain trees. We highlight results on burning in these classes in the present section.

\subsection{Graph products}
Graph products form new graphs from existing ones, and so it is natural to study how the burning number in this context. We first recall several well-known graph products. Let $G$ and $H$ be graphs, which are called \emph{factors}. Define the \emph{Cartesian product} of $G$ and $H$, written $G\square H,$ to have vertices $V(G)\times V(H),$ and vertices $(a,b)$ and $(c,d)$
are adjacent if $a=c$ and $bd \in E(H)$ or $ac \in E(G)$ and $b=d.$ Define the \emph{strong product} of $G$ and $H$, written $G \boxtimes H,$ to have vertices $V(G)\times V(H),$ and vertices $(a,b)$ and $(c,d)$
are adjacent if $a=c$ and $bd \in E(H),$ $ac \in E(G)$ and $b=d,$ or $ac \in E(G)$ and $bd \in E(H).$  The \emph{lexicographic product} of $G$ and $H$, written $G\circ H,$ has vertices $V(G)\times V(H),$ and vertices $(a,b)$ and $(c,d)$
are adjacent if $ac \in E(G)$  or $a=c$ and $bd \in E(H).$

The following theorem gives bounds on the burning number of Cartesian and strong products of graphs in terms of the bounding number of their factors.

\begin{theorem}[\cite{prod}]\label{tprod}
If $G$ and $H$ are connected graphs, then we have that
$$\max\{b(G),b(H)\} \le b(G \boxtimes H) \le b(G\square H) \le \min\{b(G)+\mathrm{rad}(H),b(H)+\mathrm{rad}(G)\}.$$
\end{theorem}

An important class of graph products are grids, which are products of paths. The $m \times n$ \emph{Cartesian grid}, defined as $P_m \square P_n,$ is denoted by $G_{m,n}.$ The value of $b(G_{m,n})$ for $m$ a function of $n$ was first studied in \cite{rburn}.

\begin{theorem}[\cite{rburn}] \label{thm grid burning number}
For $m = m(n)$, \[b\big(G_{m,n}\big) = \begin{cases}
(1 + o(1))\sqrt[3]{\frac{3}{2} mn} \: & \mbox{ if }n \geq m = \omega\big(\sqrt{n}\big), \\[0.5cm]
\Theta\big(\sqrt{n}\big) \: & \mbox{ if } m = O\big(\sqrt{n}\big).
\end{cases} \]
\end{theorem}
While Theorem~\ref{thm grid burning number} gives an asymptotically tight value for the burning number of grids where $n \geq m = \omega\big(\sqrt{n}\big),$ only the growth rate is given in the remaining case where $m = O\big(\sqrt{n}\big)$.
We refer to the family of grids $b(G_{c\sqrt{n},n})$ for constant $c > 0$ as \emph{fences}, as they are by definition wider than they are tall.

The following theorem from \cite{bonf} improves on Theorem~\ref{thm grid burning number}, giving explicit lower and upper bounds on the burning number of fences.
\begin{theorem}[\cite{bonf}] \label{thm main result}
Let $c > 0$. If $\ell = \max\{k\in\mathbb{N} : (k - 1)\sqrt{kn} + 1 \leq c\sqrt{n}\}$, then we have that
\[
b(c\sqrt{n},n) \geq \begin{cases}
\big(1 + o(1) \big) \left(\frac{c}2+\sqrt{1-\frac{c^2}4}\right) \sqrt{n} \,, & \mbox{ if } 0 < c < 2, \\[0.5cm]
\sqrt{\ell n} \,, & \mbox{ if } c \geq 2.
\end{cases}
\]
If $\ell=\lceil (c/2)^{2/3}\rceil,$ then we have that
\[
b(G_{c\sqrt{n},n})\leq 2\sqrt{\ell n}+\ell-1,
\]
and for $0<c\leq 2\sqrt{2}$, we have that
\[
b(G_{c\sqrt{n},n})\leq (1+o(1))\left(\frac{c}{2}+\sqrt{1-\frac{c^2}{16}}\right)\sqrt{n}.
\]
\end{theorem}
Another well-known graph formed from the Cartesian product are hypercubes. The $n$-dimensional \emph{hypercube}, written $Q_n,$ is the iterated Cartesian product of $n$-copies of $K_2$. In \cite{prod}, it was shown that $b(Q_n) \sim n/2.$

A \emph{strong grid} is a strong product of paths. For strong grids, we have the following asymptotic results.
\begin{theorem}[\cite{prod}] \label{tstrong}
\[
b\big(P_m\boxtimes P_n \big) = \begin{cases}
\big(1 + o(1) \big)\sqrt[3]{\frac{3}{4} mn} \: & \mbox{ if } m = \omega\big(\sqrt{n}\big), \\[0.5cm]
\Theta\big(\sqrt{n}\big) \: & \mbox{ if } m = O\big(\sqrt{n}\big).
\end{cases}
\]
\end{theorem}
As in the case for Cartesian grids, only the growth rate is provided here if $m = O\big(\sqrt{n}\big)$. Finding explicit bounds for such \emph{strong fences} has not yet been investigated.

For lexicographic products, we have the following bounds. Note that in the case that $G$ is a single vertex, then $b(G\circ H) =b(H).$
\begin{theorem}[\cite{prod}] \label{tlex}
Let $G$ be a connected graph with order at least two and $H$ any graph. We then have that
$$b(G) \le b(G \circ H) \le b(G) + 1.$$
\end{theorem}

\subsection{Binomial Random Graphs}
Let $0 \le p \le 1$ and let $\Omega$ be the family of all graphs on $n$ vertices. To every graph $G \in \Omega$ we assign a probability
$$
\Prob ( \{G\} ) = p^{|E(G)|} (1-p)^{ {n \choose 2} - |E(G)|}.$$
We denote this probability space by $\mathbb{G}(n,p).$ The space $\mathbb{G}(n,p)$ is often referred to as the \emph{binomial random graph}. Note also that this probability space can informally be viewed as a result of ${n \choose 2}$ independent coin flips, one
for each pair of vertices $u,v$, where the probability of adding an edge $uv$ is equal to $p$. For background on random graphs, see the books \cite{bol,BookFK,JLR}.

The burning number becomes a random variable on $\G(n,p)$. Results for the burning number of $\G(n,p)$ were provided in \cite{rburn}, and are summarized in the following theorem. We
say that an event in a probability space holds \emph{asymptotically almost surely} (\emph{a.a.s.}) if its probability tends to one as $n$ goes to
infinity.
\begin{theorem}[\cite{rburn}]\label{CH9thm:mainRandom}
Let $G \in \G(n,p)$, $\eps >0$, and $\omega = \omega(n) \to \infty$ as $n \to \infty$ but $\omega = o(\log \log n)$. Suppose first that
$$
d=d(n)=(n-1)p \gg \log n
$$
and
$$
p \le 1-(\log n + \log \log n + \omega)/n.
$$
Let $i \ge 2$ be the smallest integer such that
$$
d^i/n - 2 \log n \to \infty.
$$
The following holds a.a.s.
\begin{equation*}\label{CH9eq:i}
b(G) =
\begin{cases}
i & \text{ if } \ \ d^{i-1}/n \ge (1+\eps) \log n\\
i \text{ or } i+1 & \text{ if } \ \ \scriptstyle{(1-\eps) \log d \le d^{i-1}/n < (1+\eps) \log n}\\
i+1 & \text{ if } \ \ d^{i-1}/n < (1-\eps) \log d.
\end{cases}
\end{equation*}

If
$$
1-(\log n + \log \log n + \omega)/n < p \le 1-(\log n + \log \log n - \omega)/n,
$$
then a.a.s.\ $b(G) = 2 \text{ or } 3.$
Finally, if
$$
p > 1-(\log n + \log \log n - \omega)/n,
$$
then a.a.s.\ $b(G) = 2.$
\end{theorem}

\subsection{Other graph classes}

A \emph{path-forest} is a disjoint union of paths. If a path-forest $G$ of order $n$ has $t$ components, then observe that $b(G) \ge \max \{\lceil n^{1/2} \rceil ,t \}.$ Upper bounds on the burning number of path forests were given in \cite{bl}.

\begin{theorem}[\cite{bl}] \label{thm:path-forest}
Let $G$ be a path-forest of order $n$ with $t \ge 1$ components. We then have the following bounds:
\begin{enumerate}
\item $b(G) \le \left \lfloor \frac{n}{2t} \right \rfloor  + t.$
\item If $t \le \lceil n^{1/2} \rceil$, then $b(G) \le \left \lceil n^{1/2} + \frac {t-1} 2 \right \rceil.$
\end{enumerate}
\end{theorem}
Note that (2) improves on the bound (1) for smaller values of $t.$ Results on burning numbers of path forests were also considered in \cite{das,sim}.

The \emph{circulant graph} on $n$ vertices with \emph{distance set} $S$ has vertex set $\mathbb{Z}_n$ and edge set $\{ xy: x - y \in S\},$ where $S \subseteq \mathbb{Z}_n$ and $x \in S$ implies $-x \in S,$ with addition taken modulo $n.$
In \cite{shannon}, exact values of the burning numbers of $3$-regular circulants were found, along with bounds on the burning numbers of $4$-regular circulants.

Let $n \ge 3$ and $k$ be integers satisfying $1 \le k \le n-1.$ The \emph{generalized Petersen graph} $P(n,k)$ has vertices $\{ u_i , v_i: i=0,1,\ldots n-1 \}$ and edges (with subscripts modulo $n$) given by $\{u_i u_{i+1}, u_i v_i, v_i v_{i+k}, i=0,1, \ldots n-1\}.$ In \cite{sim}, it was proven that
$$ \ceil*{\sqrt{\floor*{\frac{n}{k}}}} \le b(P(n,k)) \le \ceil*{\sqrt{\floor*{\frac{n}{k}}}} + \floor*{\frac{k}{2}} + 2.$$

For positive integers $a,b,c$ define the \emph{theta graph} $\Theta_{a,b,c} (u,v)$ to be the graph consisting of a pair of vertices $\{ u,v \}$ and three internally-disjoint paths between them of lengths $a+1,$ $b+1,$ and $c+1.$ Note that $\Theta_{a,b,c} (u,v)$ has order $a+b+c+2$. It was proven in \cite{liu2} that theta graphs of order $q^2+r$, where $1\le r \le 2q+1,$ satisfy $q \le b(\Theta_{a,b,c} (u,v)) \le q+1.$  More detailed results for various values of parameters may be found in \cite{liu2}.

\section{Complexity of graph burning}\label{seccomp}

We provide a summary of what is known regarding the computational complexity of the graph burning decision problem. For background on complexity theory the reader is directed to \cite{sipser}, to \cite{app} for approximation algorithms, and \cite{zoo} for a comprehensive list of complexity classes.

We formalize the graph burning decision problem as follows.

\medskip

\noindent{\bf Problem: Graph Burning}

\noindent{\bf Instance:} A graph $G$ of order $n$ and an integer $k\geq 2$.

\noindent{\bf Question:} Is $b(G) \leq k$? In particular, does $G$ contain a burning sequence $(x_1, x_2, \ldots, x_k)$?

\medskip

While it will not come as a surprise that the Graph Burning problem is \textbf{NP}-complete, it is interesting that it remains so for fairly restrictive graph classes. The Graph Burning problem was shown to be \textbf{NP}-complete when restricted to trees of maximum degree three in \cite{bjr3}. Further, it is \textbf{NP}-complete when restricted to spider graphs, and also for disconnected graphs such as path-forests. In \cite{bjr3}, a polynomial time algorithm was given for finding the burning number of path-forests and spider graphs, when the number of arms and components is fixed.

In \cite{liu1} and independently in \cite{hil}, it was shown that the Graph Burning problem is \textbf{NP}-complete for caterpillars of maximum degree 3. In \cite{gupta}, it was shown that the Graph Burning problem is \textbf{NP}-complete when restricted to any one of the classes of interval graphs, permutation graphs, or disk graphs. Burning was considered for directed graphs in \cite{rj}, where it was proved that computing the burning number of a directed tree is \textbf{NP}-hard. Further, the Graph Burning problem is \textbf{W[2]}-complete for directed acyclic graphs. In \cite{kare}, it was shown that the Graph Burning problem can be solved in polynomial time on cographs and split graphs.

In \cite{bjr3}, a polynomial time approximation algorithm with approximation factor $3$ was given for general graphs. In \cite{bl}, a $\frac{3}{2}$-approximation algorithm was given for burning path-forests. In \cite{bk}, a polynomial time approximation algorithm with approximation factor $2$ was given for trees. In case the graph is a path-forest with a constant number of paths, the results of \cite{bk} provide a dynamic programming algorithm that creates an optimal solution in polynomial time. When the number of paths is not a constant, they provided two approximation schemes. The first scheme works under a regularity condition which implies the lengths of paths are asymptotically equal. For this scheme, they reduced the problem to the bin covering problem to achieve a fully polynomial time approximation scheme for the problem. For the general setting, when there is no assumption on the length of the paths, they found a polynomial time approximation scheme which runs in time polynomial in the size of the graph.

A graph decision problem is in \textbf{APX} if it is in \textbf{NP} and allows a polynomial time approximation algorithm with approximation ratio bounded by a constant. A graph decision problem is \textbf{APX}-hard if there is a polynomial time approximation scheme reduction from every problem in \textbf{APX} to that problem. In \cite{mondal}, it was proven that the Graph Burning problem is \textbf{APX}-hard, answering a question from \cite{bk}. It was also proven in \cite{mondal} that even if the burning sources are given as an input, computing a burning sequence itself is \textbf{NP}-hard.

In light of the approximation algorithms known for graph burning, heuristics for the problem were considered in \cite{g,kare}. For example, in \cite{g}, the authors introduce three heuristics based on eigenvalue centrality for graph burning: Backbone Based Greedy Heuristic, Improved Cutting
Corners Heuristic and Component Based Recursive Heuristic.

In \cite{kare}, the parameterized complexity of graph burning was studied. Several problems related to parameterized complexity of graph burning from \cite{kare} were solved in \cite{ko}. For example, they proved in \cite{ko} that the Graph Burning problem parameterized by $k$ is \textbf{W[2]}-complete.

\section{Future directions in graph burning}\label{secopen}

We finish by discussing conjectures and directions on graph burning. The main conjecture in the area is the burning number conjecture, which states that every connected graph $G$ of order $n$, $b(G)\leq \lceil {n}^{1/2}\rceil.$ While the full conjecture remains unresolved, it would be interesting to consider other classes where it holds, such as prescribed classes of trees.

An interesting direction would be to find connections with the burning number and other graph parameters. As stated earlier, an observation in \cite{BJR} was that for any graph $G$ with radius $r$ and diameter $d$, we have that
$\lceil (d+1)^{1/2}\rceil\leq b(G)\leq r+1.$ Bounds on the burning number were also provided in \cite{BJR} utilizing the $k$-distance domination number. In recent work \cite{hog}, the burning number provides bounds on the graph throttling number, which is a graph parameter related to the cop number in the game of Cops and Robbers.

On the complexity side, a conjecture from \cite{bjr3} is that for a tree $T$ of radius $r,$ we can recognize in polynomial time whether or not $b(T) = r + 1.$ In \cite{kam}, an algorithm with approximation factor of $2 + o(1)$ was provided for burning graphs of
bounded tree-length. An open problem from \cite{kam} is to find algorithms with similar performance for graph classes such as planar graphs.

We may consider variants of the graph burning process. In \cite{dm}, the concepts of fast and slow burning were introduced. In $k$-\emph{fast burning}, given a graph $G$ and $k \in \mathbb{N},$ burned vertices of $G$ spread to all their
$k$-neighbors. This reduces to ordinary graph burning when $k = 1.$ In $k$-\emph{slow burning}, in each round, burning spreads to $k$ neighbors of our choosing. It would be interesting to these variants of burning for graph classes such as trees and hypercubes. A variant of potential interest referenced in \cite{dm,mondal} would
be \emph{edge burning}, where sources are edges and spread to incident vertices or edges.

Another variant is random graph burning, where sources are chosen via a prescribed stochastic process; see \cite{rburn}. For example, we may consider a uniform choice of sources: at round $i$ of
the process, a source is selected uniformly at random from $V$ (with replacement). Let $b_R(G)$ be the random variable associated with the first round all vertices of $G$ are burned.  In \cite{rburn}, it was proven that a.a.s.\ $b_R(P_n) \sim  \sqrt{n \log n / 2}.$
The analysis of $b_R$ is open for other classes of graphs such as hypercubes. We may also consider non-uniform random processes for choosing sources such as preferential attachment, where higher degree vertices are more likely chosen as sources.

A final direction we discuss is to consider graph burning in infinite graphs. The burning number of an infinite graph may be an infinite cardinal. However, we may consider more finitary-type questions related to density considerations. This approach was taken in \cite{bgs}, where densities of burned vertices were considered in infinite Cartesian grids. In that approach, we consider growing grids in the Cartesian plane, centered at the origin.  If the grids are of height and width $2cn+1$ at round $n$, then it was shown in \cite{bgs} that all values in the real number interval $\left [ \frac{1}{2c^2} , 1 \right ]$ are possible densities for the burned set.  It would be interesting to consider this density approach for burning other infinite grids, such as strong or hexagonal grids.

\end{document}